\documentclass{article}

\def \J {\mathbf {J}}

\def \B {{\cal B}}

\def\uu{\bigsqcup}
\def\eps{\varepsilon}

\def \B {{\mathcal {B}}}
\usepackage[T1]{fontenc} 
\usepackage[cp1251]{inputenc}  
\usepackage[russian]{babel}

\title{
Типичные расширения рекуррентны  для апериодического автоморфизма вероятностного пространства }
\author{ Д.К. Козырева, В.В. Рыжиков }
\date{}

\usepackage[T1]{fontenc} 
\usepackage[cp1251]{inputenc}  
\usepackage[russian]{babel}

\usepackage{graphicx}
\graphicspath{{pictures/}}
\DeclareGraphicsExtensions{.pdf,.png,.jpg}

\begin{document}

\date{}

\maketitle

\Large
\begin{abstract}  Доказана типичность свойства рекуррентности для расширений  апериодического автоморфизма вероятностного пространства Лебега.
\vspace{2mm}

Ключевые слова: \it Автоморфизм пространства с мерой, косое произведениние, категории Бэра, типичное расширение. \rm
 \end{abstract}

\Large

\section{Введение }

Группа автоморфизмов $Aut(\mu)$ вероятностного пространства $(X, \B, \mu)$ оснащается  полной метрикой Халмоша $\rho$, определенной формулой
$$ \rho(S, T) = \sum_i 2^{-i}\ (\mu(SA_i \Delta TA_i) + \mu(S^{-1}A_i \Delta T^{-1}A_i)), $$
где $S,T \in Aut(\mu)$, а $\{A_i\}$ -- некоторое фиксированное семейство, плотное в алгебре $\B$.  
В связи с этим возникает возможность изучать типичные (в смысле Бэра) свойства автоморфизмов. Так, например, слабое перемешивание является типичным свойством, а перемешивающие  автоморфизмы образуют множество первой категории, см. \cite{H}.  Метод категорий позволил решить некоторые важные задачи спектральной теории динамических систем, см. \cite{A}, \cite{T}.

В группе автоморфизмов $Aut(\mu\times\mu)$ можно выделить замкнутое подмножество и ставить задачу о типичности   в рамках этого ограничения.  Далее нас интересует  следующий случай. Через $Ext(S)$ обозначим все косые произведения $R$  над $S$. Напомним, что такие   $R$ определены  формулой
$$R(x,y)= (Sx, T_x y), \ x\in X, \ y\in Y.$$
Подразумевается, что  $\{T_x\}$ -- измеримое семейство автоморфизмов пространства $(Y,\mu)$. Для таких используем $R$ обозначение $(S,T_x)$.  Мы рассматриваем  далее только случай, когда   $Y=X$, $\mu'=\mu$. 
Метрика Халмоша на $Aut(\mu\otimes\mu)$  индуцирует полную метрику на подпространстве  $Ext(S)$. Таким образом, имеет смысл говорить типичных в  смысле Бэра расширениях автоморфизма.

Семейство, содержащее $G_\delta$-множество, плотное в $Ext(S)$, называется типичным. Свойство расширения типично, если им обладают представители некоторого типичного множества. 
Для расширений $R=(S,T_x)$ и соответствующих коциклов
$C(x,n,R)=T_{S^{n-1}}...T_{Sx} T_x$ определим   класс  
$$Rec(S)=\{R \in Ext(S) : \forall m,N \,\, \exists n>N  \,\, \mu(D(m, n, R))>0\},$$ 
где $D(m, n, R)=\{x: \rho(C(x,n,R), Id)<1/m\}$.
Расширения из $Rec(S)$ называются рекурентными. 
Для рекуррентных $R$    для сколь угодно больших $n$ на множестве $x$-ов положительной меры коцикл $C(x,n,R)$ оказывается близок к тождественному преобпазованию.  Это можно интерпретировать как эффект возвращения к нейтральному элементу при случайном блуждании по группе $Aut(\mu)$. 
В \cite{R1} показано, что свойство рекуррентности типично для  расширений слабо перемешивающего автоморфизма. Данная статья представляет собой обобщение этого результата для  расширений апериодического и тем самым для  расширений эргодического автоморфизма. 

Отметим, что в типичном случае возвращение к нейтральному элементу, как правило, происходит на маленьком множестве $x$-ов. При этом  на большом множестве $x$-ов  автоморфизмы $C(x,n,R)$ "уходят на бесконечность".  Последнее надо понимать как асимптотическое нахождение $C(x,n,R)$ вне всякого компакта в группе $Aut(\mu)$. Сказанное вытекает из типичности свойства относительного слабого перемешивания \cite{GW} для расширений эргодических автоморфизмов.
Рекуррентность типичных коциклов существенно использовалось для поднятия свойства перемешивания \cite{R1} и кратного перемешивания вместе с  инвариантом, формулируемым в терминах самоприсоединений, см. \cite{R2}. 
При положительной энтропии эргодического базового автоморфизма типичное расширение оказывается изоморфным этому автоморфизму \cite{AW} .

\bf Теорема 1.  \it Если $S$ -- апериодический автоморфизм,  класс $Rec(S)$ содержит плотное $G_\delta$-множество. \rm

\section{Плотность класса когомологичности тривиального расширения}
 Обозначим через $\J$ семейство всех автоморфизмов пространства $(X\times Y, \mu\otimes\mu')$, 
оставляющих инвариантными  множества вида $A\times Y$ для всех $A\subset X$. Такие автоморфизмы
являются косыми произведениями над  тождественным преобразованием $Id$. 
Косые произведения $R_1, R_2$ (и соответствующие коциклы)  над автоморфизмом $S$ называются когомологичными,  если выполняется 
$R_1=JR_2J^{-1},$ для некоторого $J\in \J.$

\bf Лемма 1. \it Для апериодического автоморфизма $S$ класс когомологичности  косого произведения $S \times Id$ плотен в $Ext(S)$. \rm

 Доказательство. Фиксируем  косое произведение $R=(S,T_x)$ и $\delta>0$. Найдем отображение $\Phi\in Aut$ такое, что 
$$ \mu(x:\Phi^{-1} R_0 \Phi x \neq Rx)<\delta. $$
Воспользуемся леммой Рохлина-Халмоша:

 \it для апериодического автоморфизма $S$, числа $\eps >0$ и натурального числа $N$ существует такое измеримое  $B$, что множества $T^iB$ не пересекаются при $0 \leq i \leq N-1$, причем \rm
$$ \mu\left(\bigsqcup\limits_{i=0}^{N-1}S^iB \right)>1-\eps. $$\rm

Заметим что для множества $A=B \times Y$ верно $R^iA=R^i_0A$ для всех $i=0,...,n-1$. Выбираем  сохраняющее меру отображение  $\Phi$  такое, что  $\Phi R^ix:=R_0^i \Phi x$ для $i=0,...,n-1$. Тогда $\Phi$ определено   на всей башне, кроме оставшенося множества.  На нем  $\Phi$  доопределяется до автоморфизма  произвольным образом.


Для всех точек $x$, не принадлежащих последнему этажу и остатку, выполнено  $Rx=\Phi^{-1} R_0 \Phi x$.  Значит, $\mu(x:\Phi^{-1} R_0 \Phi x \neq Rx)<\eps+\frac{1}{m}$. Выбрав $\eps$ и $m$ такие, что  $\eps+1/m <\delta$,  мы завершаем доказательство леммы.

\section{Типичность рекуррентности}

\bf Лемма 2. \it Рекуррентнтные коциклы образуют $G_\delta$-множество.\rm

Доказательство.  Рассмотрим семейство расширений
$$C_m(S)=\{R \in Ext(S): \forall N \in \mathbf{N} \,\,\, \exists \, n>N : \mu \left( x : \rho(C(x,n,R), Id) < {1}/{m} \right) >0 \}, $$
где $\rho$ -- метрика Халмоша в $Aut(\mu)$. Напомним, что $R=(S,T_x)$, $C(x,n,R)=T_{S^{n-1}x}\dots T_{Sx}T_x$. 
Множество 
$$U_{m,n}=\{R \in Ext(S): \varphi(R)=\mu \left( x : \rho(C(x,n,R), Id) < {1}/{m} \right) >0 \}$$ при фиксированных $m,n$
 является открытым, так $\varphi(R)$ непрерывно зависит от $R$.
Таким образом, 
$$Rec=\bigcap_m C_m(S)=\bigcap_m  \bigcap_N\bigcup_{n>N} U_{m,n}$$
является $G_\delta$-множеством.

\vspace{2mm}
\bf Лемма 3. \it Класс $Rec(S)$  плотен в $Ext(S)$.\rm

\vspace{2mm}
Доказательство.  Все расширения автоморфизма $S$, когомологичные произведению $S \times Id$, образуют плотный в $Ext(S)$ класс (лемма 1). Они имеют вид $J^{-1}(S \times Id)J$,  $J\J$. В классе $Ext(S)$ в силу его сепарабельности  найдется счётное всюду плотное в $Ext(S)$  множество расширений $R_i$, имеющих вид
$$R_i=J_i^{-1}(S \times Id)J_i, \,\,\, J_i=(Id, J_{i,x}). $$
При этом можно считать,  что  для каждого $i$ семейство автоморфизмов  $\{J_{i,x}: x\in X\}$  счетно.  Здесь $X=\uu_k B_{i,k}$, и для каждого $k$  автоморфизмы $J_{i,x}$ одинаковы на измеримом множестве $B_{i,k}$. Назовем 
такие семейства  простыми. Сказанное является известным фактом, который использовался в \cite{R1}. Ограничимся коротким пояснением того, как построить такое семейство.  
В силу сепарабельности $Ext(S)$    множество расширений автоморфизма $S$, когомологичные произведению $S \times Id$, тоже сепарбельно. 
Возьмем в нем счетное всюду плотное множество. Каждое расширение из этого множества приближается последовательностью простых  расширений. Объединим все эти последовательности в требуемое счетное семейство $\{R_i\}$. 

Покажем, что $R_i\in Rec(S)$. Рассмотрим множество $B_{i,k}$ положительной
меры. Для всякого $N$ найдется $n>N$ такое, что $\mu (S^nB_{i,k}\cap B_{i,k})>0$. 
Для  $x\in \ B_{i,k}\cap S^{-n}B_{i,k}$ непосредственно проверяется, что $C(x,n,R_i)=Id$.  А это означает, что $R_i\in Rec(S)$.
Лемма доказана. С учетом леммы 2 получаем утверждение  теоремы 1.

На самом деле мы доказали чуть больше: типичное расширение рекуррентно  на произвольном  множестве $A\subset X$ положительной меры. 
Пусть $$D(m, n, R, A)=\{x\in A: \rho(C(x,n,R), Id)<1/m\},$$
где, напомним,  $C(x,n,R)=T_{S^{n-1}}...T_{Sx} T_x$, $R=(S,T_x)$.
Обозначим 
$$Rec(S, A)=\{R \in Ext(S) : \forall m,N \,\, \exists n>N  \,\, \mu(D(m, n, R,A))>0\},$$ 
где $D(m, n, R, A)=\{x\in A: \rho(C(x,n,R), Id)<1/m\}$.
Расширения из $Rec(S, A)$ называются рекурентными на $A$.  

\vspace{2mm}
\bf Теорема 2. \it
Для апериодического  преобразования  $S$ и множестве $A$ положительной меры
класс рекурентных на $A$  расширений автоморфизма $S$  содержит плотное $G_\delta$-множество. \rm

\vspace{2mm}
Доказательство повторяет изложенные выше аргументы с той лишь разницей, что 
в аналоге леммы 3 надо сказать,  что рассматриваются  множества $B_{i,k}$, лежащие в $A$.

\end{document}